# ADDENDUM TO THE DISCUSSION OF "BREAKDOWN AND GROUPS"[1]

BY P. LAURIE DAVIES AND URSULA GATHER

*University of Duisburg–Essen and Technical University Eindhoven, and University of Dortmund*

In his discussion of Davies and Gather [*Ann. Statist.* **33** (2005) 977–1035] Tyler pointed out that the theory developed there could not be applied to the case of directional data. He related the breakdown of directional functionals to the problem of definability. In this addendum we provide a concept of breakdown defined in terms of definability and not in terms of bias. If a group of finite order $k$ acts on the sample space we show that the breakdown point can be bounded above by $(k-1)/k$. In the case of directional data there is a group of order $k=2$ giving an upper bound of $1/2$.

**1. Introduction.** It has been argued that breakdown occurs in directional data (see [4]) when contamination causes the direction to change by $180°$ (see, e.g., [2, 3]). More formally this definition of breakdown point of a directional functional $T$ at a distribution $P$ may be written as

$$(1) \quad \begin{aligned} \varepsilon^*(T, P, d) \\ = \inf\{\varepsilon > 0 : |T(P) - T(Q)| = \pi \text{ for some } Q \text{ with } d(P, Q) < \varepsilon\}. \end{aligned}$$

Similarly in the case of linear spaces breakdown can be said to occur when the linear space derived from the contaminated data is orthogonal to that derived from the uncontaminated data. As noted by Tyler in the discussion this definition of breakdown cannot be formulated within the theory of Davies and Gather [1], which requires unbounded metrics deriving from certain banned or forbidden parameter values such as 0 and $\infty$ in the scale case. In this context breakdown is said to occur when the realized parameter values tend to or take on a banned value. For directional data there will, in

Received September 2004; revised September 2005.
[1]Supported in part by Sonderforschungsbereich 475, University of Dortmund.
*AMS 2000 subject classifications.* Primary 62G35; secondary 62H11.
*Key words and phrases.* Breakdown, equivariance, directional data.

general, be no banned direction so there can be no concept of breakdown based on banned values. Nevertheless if the data is such that small changes can result in a change of direction of 180° the directional functional is highly nonstable. As pointed out by Tyler in the discussion this behavior is related to the problem of definability of the functional. Tyler writes "*This implies that breakdown occurs when a 'well-defined' vector becomes 'undefined'*". We take up this idea and use it to define what we call the "definability breakdown point".

## 2. Breakdown and invariant distributions.

2.1. *Definition of breakdown.* From Theorem 3 of [3] and Theorem 4.1 of [2] it follows that the breakdown point (1) of the circular mean at a distribution $P$ is at most $1/2$ and can be arbitrarily small. He and Simpson [2] write "the breakdown point must go to 0 as the data become more dispersed over the sphere" and "it is possible to have a breakdown point of near $1/2$ if the data are concentrated." We give a definition of breakdown point which makes no use of bias but helps to understand why, for rotation equivariant functionals, a high breakdown point can be attained at concentrated distributions but only low ones are possible at disperse distributions. Furthermore it explains why $1/2$ is the largest possible breakdown point for such functionals.

We use the notation of Davies and Gather [1]. A functional $T:\mathcal{P}_T \to \Theta$, $\mathcal{P}_T \subset \mathcal{P}$ is called equivariant if the following hold:

(a) $\mathcal{P}_T$ is closed under all $g \in \mathcal{G}$,
(b) $T(P^g) = h_g(T(P))$ for all $P \in \mathcal{P}_T$ and $g \in \mathcal{G}$,
(c) $T$ is well defined on $\mathcal{P}_T$.

This leads to the following definition of the definability breakdown point:

$$(2) \qquad \varepsilon^*(T, P, d) = \inf\{\varepsilon > 0 : d(P, Q) < \varepsilon \text{ for some } Q \notin \mathcal{P}_T\},$$

with of course $\varepsilon^*(T, P, d) = 0$ if $P \notin \mathcal{P}_T$.

2.2. *Invariant measures.* We call a measure $P$ invariant if $P^g = P$ for some $g$ with $h_g \neq h_\iota$. We denote the set of all such distributions by $\mathcal{P}_{\text{inv}}$. If $T$ is equivariant and $P \in \mathcal{P}_{\text{inv}}$ we have $T(P) = T(P^g) = h_g(T(P))$, which is not possible as $h_g \neq h_\iota$. This implies

$$(3) \qquad \mathcal{P}_{\text{inv}} \subset \mathcal{P} \setminus \mathcal{P}_T$$

for every equivariant functional $T$ and hence

$$(4) \qquad \varepsilon^*(T, P, d) \leq \inf\{\varepsilon > 0 : d(P, Q) < \varepsilon \text{ for some } Q \in \mathcal{P}_{\text{inv}}\}.$$

In the case of directional data the uniform distribution $U$ belongs to $\mathcal{P} \setminus \mathcal{P}_T$ for any equivariant functional $T$ and hence we have the upper bound

$$\varepsilon^*(T, P, d) \leq d(P, U). \tag{5}$$

This implies that the breakdown point is low at very disperse distributions, that is, those close to the uniform distribution.

2.3. *Finite sub-groups.* Suppose $\mathcal{G}$ contains a finite sub-group $\mathcal{G}_k$ of order $k \geq 2$ so that $g^k = \iota$ for all $g \in \mathcal{G}_k$. For any distribution $P$ we set

$$P_k = \frac{1}{k} \sum_{j=0}^{k-1} P^{g^j}. \tag{6}$$

Then $P_k^g = P_k$ so that $P_k \in \mathcal{P}_{\text{inv}}$ and hence

$$\varepsilon^*(T, P, d) \leq d(P, P_k).$$

If the metric $d$ satisfies (2.1) and (2.2) of [1] we have

$$d(P, P_k) \leq \frac{k-1}{k},$$

where

$$\tilde{P}_k = \frac{1}{k-1} \sum_{j=1}^{k-1} P^{g^j}.$$

Examples of sample spaces and groups $\mathcal{G}$ with subgroups of order $k = 2$ are the unit circle and the unit sphere in any dimension. In all cases the maximum definability breakdown point of any direction functional at a measure $P$ is at most $1/2$ and may be much smaller as implied by (5).

2.4. *Total variation and finite sample breakdown points.* Although the above results can be extended to finite sample breakdown points this may not always make sense. In particular if $P$ is an empirical measure then the breakdown point measured using the total variation metric $d = d_{\text{tv}}$ may be reduced from $1/2$ to $1/n$ by the smallest of alterations in the values of the data. The same applies to the finite sample breakdown points. This is the only example we know where the use of a metric which allows for minor alterations in the values of the data points leads to a completely different breakdown point.

**Acknowledgment.** We wish to acknowledge the remarks of an anonymous referee whose comments lead to a considerable increase in clarity.

FACHBEREICH 06—MATHEMATIK
  UND INFORMATIK
UNIVERSITÄT DUISBURG–ESSEN
45117 ESSEN
GERMANY
E-MAIL: laurie.davies@uni-duisburg-essen.de

FACHBEREICH STATISTIK
UNIVERSITÄT DORTMUND
44221 DORTMUND
GERMANY
E-MAIL: gather@statistik.uni-dortmund.de